\documentclass[12pt]{article}
\usepackage{amstext}
\usepackage{amsfonts}
\usepackage{amssymb}
\usepackage{amsbsy}
\usepackage{latexsym}
\usepackage{xy}
\usepackage{hhline}
\xyoption{all}
\newcounter{num} %

\newenvironment{theo}
{\refstepcounter{num}%
\bigskip\noindent{\bf Theorem~\arabic{num}. }\it}

\newenvironment{cor}
{\refstepcounter{num}%
\bigskip\noindent{\bf Corollary~\arabic{num}. }\it}

\newenvironment{lemma}
{\refstepcounter{num}%
\bigskip\noindent{\bf Lemma~\arabic{num}. }\it}

\newenvironment{remark}
{\refstepcounter{num}%
\bigskip\noindent{\bf Remark~\arabic{num}. }\it}

\newcommand{\Ref}[1]{(\ref{#1})}

\newcounter{theexample}
\newcommand{\example}
{\refstepcounter{theexample}%
\bigskip\noindent{\bf Example~\arabic{theexample}.}}

\newenvironment{proof}{\medskip\noindent{\it Proof. }}
{$\Box$ \bigskip}

\newenvironment{eq}{\begin{equation}}{\end{equation}}
\newcommand{\si}{\sigma}

\newcommand{\be}{\beta}
\newcommand{\ga}{\gamma}
\newcommand{\la}{\lambda}

\newcommand{\La}{\Lambda}

\newcommand{\ov}[1]{\overline{#1}}

\newcommand{\tr}{{\rm tr}}

\newcommand{\image}{\mathop{\rm Im}}
\newcommand{\sign}{\mathop{\rm{sgn }}}

\newcommand{\un}[1]{{\underline{#1}} }

\newcommand{\Pf}{{\mathop{\rm{pf }}}}
\renewcommand{\P}{{\mathop{\ov{\rm{pf}}}}}
\newcommand{\DP}{{\rm DP} }
\newcommand{\F}{{\mathop{\rm{bpf }}}}

\newcommand{\NN}{{\mathbb{N}} }
\newcommand{\ZZ}{{\mathbb{Z}} }
\newcommand{\QQ}{{\mathbb{Q}} }
\newcommand{\M}{{\mathcal M} }
\newcommand{\ovphi}[1]{\varphi(#1)}

\newcommand{\Lus}{\{}        
\newcommand{\Rus}{\}_{m}}  %

\newcommand{\Symmgr}{\mathcal{S}}                 

\begin{document}
\title{On block partial linearizations of the pfaffian.}
 \author{
A.A. Lopatin \\
Institute of Mathematics, \\
Siberian Branch of \\
the Russian Academy of Sciences, \\
Pevtsova street, 13,\\
Omsk 644099 Russia \\
artem\underline{ }lopatin@yahoo.com \\
http://www.iitam.omsk.net.ru/\~{}lopatin/\\
}
\date{} 
\maketitle

\begin{abstract}
Amitsur's formula, which expresses $\det(A+B)$ as a polynomial in coefficients of the
characteristic polynomial of a matrix, is generalized for partial linearizations of the
pfaffian of block matrices. As applications, in upcoming papers we determine generators
for the $SO(n)$-invariants of several matrices and relations for the $O(n)$-invariants of
several matrices over a field of arbitrary characteristic.
\end{abstract}



\section{Introduction}\label{section_intro}

Throughout, $K$ is an infinite field of arbitrary characteristic.
All vector spaces, algebras, and modules are over $K$ unless
otherwise stated.

A block partial linearization of the pfaffian (b.p.l.p.) is a
partial linearization of the pfaffian of block matrices. The
concept of b.p.l.p.~is important for the invariant theory.
Examples of b.p.l.p.~include the function $\DP$, of three
matrices, used to describe generators for semi-invariants of mixed
representations of quivers (see~\cite{LZ1}), the matrix function
$\si_{t,s}$, which describes relations between generators for the
invariants of mixed representations of quivers
(see~\cite{ZubkovI}), and a partial linearization of the
determinant of block matrices as well as the determinant
itself~(see part~2 of Example~\ref{ex1} for details).

Denote coefficients in the characteristic polynomial
of an $n\times n$ matrix $X$ by $\si_k(X)$, i.e., %
$$\det(\la E-X)=\la^n-\si_1(X)\la^{n-1}+\cdots+(-1)^n\si_n(X).$$ %
Recall that Amitsur's formula from~\cite{Amitsur} expresses $\si_k(A+B)$ as a polynomial
in $\si_i(a)$, where $1\leq i\leq k$ and $a$ is a monomial in $A$ and $B$ (see
Corollary~\ref{cor_Amitsur} for an explicit formulation). The main result of this paper
is the decomposition formula that generalizes Amitsur's formula for a b.p.l.p.~and that
links a b.p.l.p.~with coefficients in the characteristic polynomial of a matrix. This
answers the question formulated at the end of Section~$4$ in~\cite{ZubkovII} as a key
consideration for the following problem.

Let $K_{n,d}=K[x_{ij}(r)\,|\,1\leq i,j\leq n,\,1\leq r\leq d]$ be a polynomial
$K$-algebra and let $X_{r}=(x_{ij}(r))_{1\leq i,j\leq n}$ be an $n\times n$ matrix. Note
that the subalgebra generated by $\si_k(X_{r_1}\cdots X_{r_s})$ for $1\leq k\leq n$ and
$1\leq r_1,\ldots,r_s\leq d$ is the algebra of invariants of several $n\times n$ matrices
(see~\cite{Donkin92a}), and the subalgebra generated by $\si_k(Z_{r_1}\cdots Z_{r_s})$,
where $Z_r$ is $X_r$ or its transpose $X_r^t$ and $r_1,\ldots,r_s$ are in the same range,
is the algebra of $O(n)$-invariants of several matrices (see~\cite{Zubkov99}). Amitsur's
formula plays a crucial role for the description of relations between
$\si_k(X_{r_1}\cdots X_{r_s})$ in~\cite{Zubkov96}. The decomposition formula plays the
same role for the elements $\si_k(Z_{r_1}\cdots Z_{r_s})$ as it was pointed out
in~\cite{ZubkovII}. As a consequence of the decomposition formula we will describe
relations for the $O(n)$-invariants when the characteristic of $K$ is different from $2$
in the upcoming paper~\cite{Lop_relations}.

The decomposition formula will be also used to obtain generators for the
$SO(n)$-invariants of several matrices when the characteristic of $K$ is different from
$2$ (see~\cite{Lop_so_inv}). Hence we will complete description of generators for the
invariants of several matrices, that was originated by Procesi in 1976
(see~\cite{Procesi76}). Moreover, generators for the invariants under the action of a
product of classical groups on the space of (mixed) representations of a quiver will be
also established in~\cite{Lop_so_inv}.


An immediate consequence of the decomposition formula is
Corollary~\ref{cor_PP} which states that the product of two
pfaffians is a polynomial in $\si_k$.

The paper is organized as follows.

Section~\ref{section_prelim} contains general notation.

Section~\ref{section_def} starts with basic definitions of the pfaffian and a b.p.l.p.
The key notion of a tableau with substitution is introduced in order to provide a more
convenient language for dealing with b.p.l.p.-s (see Lemma~\ref{lemma_bplp}). All results
are formulated using the notion of a tableau with substitution. Section~\ref{section_def}
terminates with the short version of the decomposition formula (see
Theorem~\ref{theo_short}).

The decomposition formula (Theorem~\ref{theo_decomp}) and all
necessary definitions are stated in
Section~\ref{section_decomp_form}.

The proof of the formula spans the next three sections.

Some applications of the decomposition formula are presented in
Section~\ref{section_examples}.

\section{Preliminaries}\label{section_prelim}

In what follows, all matrices have entries in some commutative unitary algebra over the
field $K$. Denote the set of non-negative integers by ${\NN}$, the set of integers by
$\ZZ$, and by $\QQ$ the quotient field of the ring $\ZZ$.

Introduce the lexicographical order $<$ on $\NN\times\NN$ by
$(a_1,b_1)<(a_2,b_2)$ if and only if $a_1<a_2$ or $a_1=a_2$ and
$b_1<b_2$. Denote the composition of functions $f$ and $g$ by
$f\circ g(x)=f(g(x))$ or just by $fg(x)$.

The cardinality of a set $S$ is denoted by $\#S$ and $S|_{a\to b}$ denotes the result of
substitution of $b$ for $a$ in $S$. A finite sequence $a_1,\ldots,a_r$ is denoted by
$(a_1,\ldots,a_r)$ or just by $a_1\cdots a_r$.  In the latter case, $a_1\cdots a_r$ is
called a word in letters $a_1,\ldots,a_r$. Denote the degree of a word $w$ in a letter
$a$ by $\deg_a w$. Similarly to the cardinality of a set, we denote the length of a
sequence $\un{\la}=(\la_1,\ldots,\la_p)$ by $\#\un{\la}=p$. We use notation
$\Lus\ldots\Rus$ for multisets, i.e., given an equivalence $=$ on a set $S$ and
$a_1,\ldots,a_p,b_1,\ldots,b_q\in S$, we write $\Lus a_1,\ldots,a_p\Rus=\Lus
b_1,\ldots,b_q\Rus$ if and only if $p=q$ and
$$\#\{1\leq j\leq p\,|\,a_j=a_i\}=\#\{1\leq j\leq p\,|\,b_j=a_i\}$$ %
for any $1\leq i\leq p$. For integers $i<j$ denote the interval
$i,i+1,\ldots,j-1,j$ by $[i,j]$.

Denote by $id_n$ the identity permutation of the group $\Symmgr_n$
of permutations on $n$ elements. Given $i\neq j$, $(i,j)$ is the
transposition interchanging $i$ and $j$.

An element $u$ of a free monoid is called {\it primitive} if
$u\neq u_0\cdots u_0$, where the number of factors is more than
$1$.

A vector $\un{\la}=(\la_1,\ldots,\la_l)\in\NN^l$ satisfying
$\la_1\geq\cdots \geq\la_l$ and $\la_1+\cdots+\la_l=t$ is called a
{\it partition} of $t$ and is denoted by $\un{\la}\vdash t$. A
{\it multi-partition} $\un{\la}\vdash \un{t}$ is a $q$-tuple of
partitions $\un{\la}=(\un{\la}_1,\ldots,\un{\la}_q)$, where
$\un{\la}_i\vdash t_i$, $\un{t}=(t_1,\ldots,t_q)\in\NN^q$.

\section{Definitions}\label{section_def}

Let us recall definitions of the pfaffian and its partial
linearizations. 

If $n$ is even, then the pfaffian of an $n\times n$
skew-symmetric matrix $X=(x_{ij})$ is %
$$\Pf(X)=\sum\sign(\pi)\prod\limits_{i=1}^{n/2} x_{\pi(2i-1),\pi(2i)},$$%
where the sum ranges over permutations $\pi\in \Symmgr_{n}$ satisfying
$\pi(1)<\pi(3)<\cdots<\pi(n-1)$ and $\pi(2i-1)<\pi(2i)$ for all $1\leq i\leq n/2$. Define
the {\it generalized pfaffian} of an arbitrary $n\times n$ matrix $X=(x_{ij})$ by
$$\P(X)=\Pf(X-X^{t}).$$
By abuse of notation we will refer to $\P$ as the pfaffian. For $K=\QQ$ there is a more
convenient formula
\begin{eq}\label{eq_P}
\P(X)=\Pf(X-X^{t})=\frac{1}{(n/2)!}\sum\limits_{\pi\in \Symmgr_{n}}\sign(\pi)
\prod\limits_{i=1}^{n/2} x_{\pi(2i-1),\pi(2i)}.%
\end{eq} %
Note that for an $n\times n$ skew-symmetric matrix
$X=(x_{ij})$ we have %
$$\P(X)=2^{n/2}\Pf(X).$$

For $n\times n$ matrices $X_1=(x_{ij}(1)),\ldots,X_s=(x_{ij}(s))$
and positive integers $k_1,\ldots,k_s$, satisfying
$k_1+\cdots+k_s=n/2$, consider the polynomial $\P(x_1
X_1+\cdots+x_s X_s)$ in the variables $x_1,\ldots,x_s$. The
partial linearization $\P_{k_1,\ldots,k_s}(X_1,\ldots,X_s)$ of the
pfaffian is the coefficient of $x_1^{k_1}\cdots x_s^{k_s}$ in this
polynomial.
In other words, for $K=\QQ$ we have %
\begin{eq}\label{eq_lin_P}
{\P}_{k_1,\ldots,k_s}(X_1,\ldots,X_s)=
\frac{1}{c}\sum\limits_{\pi\in \Symmgr_{n}}\sign(\pi)
\prod\limits_{j=1}^s\;\prod\limits_{i=k_1+\cdots+k_{j-1}+1}^{k_1+\cdots+k_{j}}
x_{\pi(2i-1),\pi(2i)}(j),%
\end{eq} %
where $c=k_1!\cdots k_s!$. The partial linearization
$\det_{k_1,\ldots,k_s}(X_1,\ldots,X_s)$ of the determinant is
defined analogously, where $X_1,\ldots,X_s$ are $n\times n$
matrices, $k_1+\cdots+k_s=n$, and $n$ is arbitrary.

Now we can give the definition of a b.p.l.p.

Fix $\un{n}=(n_1,\ldots,n_m)\in\NN^{m}$, where $n=n_1+\cdots+n_m$
is even. For any $1\leq p,q\leq m$ and an $n_p\times n_q$ matrix
$X$ denote by $X^{p,q}$ the $n\times n$ matrix, partitioned into
$m\times m$ number of blocks, where the block in the $(i,j)$-th
position is an $n_i\times n_j$ matrix; the block in the $(p,q)$-th
position is equal to $X$, and the rest of blocks are zero
matrices.

Let $1\leq p_1,\ldots,p_s,q_1,\ldots,q_s\leq m$, let $X_j$ be an $n_{p_j}\times n_{q_j}$
matrix for any $1\leq j\leq s$, and let $k_1,\ldots,k_s$ be positive integers, satisfying
$k_1+\cdots+k_s=n/2$. The element %
\begin{eq}\label{eq_block_linearization}
{\P}_{k_1,\ldots,k_s}(X_1^{p_1,q_1},\ldots,X_s^{p_s,q_s})
\end{eq}%
is a partial linearization of the pfaffian of block matrices
$X_1^{p_1,q_1},\ldots,X_s^{p_s,q_s}\!\!$, and it is  called a {\it
block partial linearization of the pfaffian} ({\it b.p.l.p.}).
We consider only b.p.l.p.-s, satisfying %
\begin{eq}\label{eq_block_condition}
\sum_{1\leq j\leq s,\;p_j=i}k_j + %
\sum_{1\leq j\leq s,\;q_j=i}k_j=n_i\;{\rm for\;all}\; 1\leq i\leq
m,
\end{eq}%
since only these elements appear in the invariant theory in the
context of Section~\ref{section_intro}.

\example\label{ex0} Let $\un{n}=(5,3)$ and let $X_1,X_2,X_3$
be $5\times 3$, $5\times 5$, and $3\times 3$ matrices respectively. Then %
$$X_1^{1,2}=\left(
\begin{array}{cc}
0&X_1\\
0&0\\
\end{array}
\right),\quad %
X_2^{1,1}=\left(
\begin{array}{cc}
X_2&0\\
0&0\\
\end{array}
\right),\quad %
X_3^{2,2}=\left(
\begin{array}{cc}
0&0\\
0&X_3\\
\end{array}
\right)
$$
are $8\times 8$ matrices. The
b.p.l.p.~$f={\P}_{1,2,1}(X_1^{1,2},X_2^{1,1},X_3^{2,2})$ satisfies
condition~\Ref{eq_block_condition}. It is convenient to display
the data that determine $f$ as a two-column tableau filled with
arrows: %
$$\begin{array}{|c|c|}
\hline a-\!\!\!\! &\!\!\!\rightarrow\,\\
\hline b\vspace{-0.2cm}&d\\
\vspace{-0.3cm}\downarrow& \downarrow\\
\hline &\\
\hline c\vspace{-0.2cm} \\
\vspace{-0.3cm}\downarrow \\
\hhline{-~} \\
\hhline{-~}
\end{array}$$%
\noindent The arrow $a$ determines the block matrix $X_1^{1,2}$ as follows: $a$ goes from
the 1st column to the 2nd column and we assign the matrix $X_1$ to $a$. To determine the
block matrix $X_2^{1,1}$ we take {\it two} arrows $b,c$ that go from the 1st column to
the 1st column, since the degree of $f$ in entries of $X_2$ is $2$. Assign the matrix
$X_2$ to $b$ and to $c$. Finally, the arrow $d$ determines the block matrix $X_2^{2,2}$:
$d$ goes from the 2nd column to the 2nd column and $X_2$ is assigned to it.
Condition~\Ref{eq_block_condition} implies that for every column the total amount of
arrows that start or terminate in the column is equal to the length of the column. Note
that this conditions do not uniquely determine arrows of a tableau.
\bigskip

Below we formulate the definition of a tableau with substitution that gives alternative
way to work with b.p.l.p.-s. Acting in the same way as in Example~\ref{ex0}, for every
b.p.l.p.~$f$ given by~\Ref{eq_block_linearization} we construct a tableau with
substitution $(T,(X_1,\ldots,X_s))$ and define the function ${\F}_T$ such that $f=\pm
{\F}_T(X_1,\ldots,X_s)$.

\bigskip
\noindent{\bf Definition (of shapes).} A {\it shape} of dimension
$\un{n}=(n_1,\ldots,n_m)\in\NN^{m}$ is a collection of $m$ columns of cells. The columns
are numbered by $1,2,\ldots,m$, and the $i$-th column contains exactly $n_i$ cells, where
$1\leq i\leq m$. Numbers $1,\ldots,n_i$ are assigned to the cells of the $i$-th column,
starting from the top. As an example, the shape of dimension
$\un{n}=(3,2,3,1,1)$ is%
$$\begin{tabular}{ccccc} %
1&2&3&4&5  \\
\end{tabular}$$ %
\vspace{-.7cm}
$$\begin{tabular}{|c|c|c|c|c|} %
\hline 1&1&1&1&1       \\
\hline 2&2&2         \\
\hhline{---~~} 3&&3   \\
\hhline{-~-~~}
\end{tabular}$$ %

\bigskip
\noindent{\bf Definition (of a tableau with substitution).} Let
$\un{n}=(n_1,\ldots,n_m)\in\NN^{m}$ and let $n=n_1+\cdots+n_m$ be even. A pair
$(T,(X_1,\ldots,X_s))$ is called a {\it tableau with substitution} of dimension $\un{n}$
if
\begin{enumerate}
\item[$\bullet$] $T$ is a shape of dimension $\un{n}$ that is filled with arrows. An {\it
arrow} goes from one cell of the shape into another one, and each cell of the shape is
either the head or the tail of one and only one arrow.  We refer to $T$ as a {\it
tableau} of dimension $\un{n}$, and we write $a\in T$ for an arrow $a$ from $T$. Given an
arrow $a\in T$, denote by $a'$ and $a''$ the columns containing the head and the tail of
$a$, respectively. Similarly, denote by $'a$ the number assigned to the cell containing
the head of $a$, and denote by $''a$ the number assigned to the cell containing the tail
of $a$. Schematically this is depicted as
$$\begin{array}{cccc} %
\;\;&a''&\! a'&  \\
\end{array}$$ %
\vspace{-.6cm}
$$\begin{array}{c|c|c|c} %
\hhline{~--~} \vspace{-0.25cm} &&&{}'a\\
              \vspace{-0.25cm} &&\hspace{-0.45cm} \nearrow&\\
\hhline{~--~}''a &a&\;\;&\\
\hhline{~--~}
\end{array}
$$%

\item[$\bullet$] $(X_1,\ldots,X_s)$ is a sequence of matrices,
where the $(p,q)$-th entry of $X_j$ is $(X_j)_{pq}$ ($1\leq j\leq
s$). We assume that an $n_{a''}\times n_{a'}$ matrix
$X_{\ovphi{a}}$ is assigned to an arrow $a\in T$, where
$1\leq\ovphi{a}\leq s$, and $\{\ovphi{a}\,|\,a\in T\}=[1,s]$.
Moreover, we assume that if $a,b\in T$ and $\ovphi{a}=\ovphi{b}$,
then $a'=b'$, $a''=b''$.
\end{enumerate}

\bigskip
\noindent{\bf Definition (of ${\F}_T(X_1,\ldots,X_s))$.} 
Let $(T,(X_1,\ldots,X_s))$ be a tableau with substitution of dimension $\un{n}$. Define
$${\F}_T(X_1,\ldots,X_s)=\sum
\sign(\pi_1)\cdots\sign(\pi_m)\prod_{a\in T}
(X_{\ovphi{a}})_{\pi_{a''}(''a),\pi_{a'}('a)},$$ %
where the sum ranges over permutations $\pi_1\in
\Symmgr_{n_1},\ldots,\pi_m\in \Symmgr_{n_m}$ such that for any
$a,b\in T$ the conditions $\ovphi{a}=\ovphi{b}$ and $''a<{}''b$
imply that $\pi_{i}({}''a)<\pi_{i}({}''b)$ for $i=a''=b''$. Note
that each coefficient in ${\F}_T(X_1,\ldots,X_s)$, considered as a
polynomial in entries of $X_1,\ldots,X_s$, is $\pm1$. For $K=\QQ$
there is a more convenient formula
$${\F}_T(X_1,\ldots,X_s)=\frac{1}{c_T}{\F}^0_T(X_1,\ldots,X_s),$$ %
where %
\begin{eq}\label{eq_def_F0}
{\F}_T^0(X_1,\ldots,X_s)=\sum_{\pi_1\in
\Symmgr_{n_1},\ldots,\pi_m\in \Symmgr_{n_m}}
\sign(\pi_1)\cdots\sign(\pi_m)\prod_{a\in T}
(X_{\ovphi{a}})_{\pi_{a''}(''a),\pi_{a'}('a)},
\end{eq}%
and the coefficient $c_T$ is equal to
$$\prod_{j=1}^s\#\{a\in T\,|\,\ovphi{a}=j\}!.$$ %

\begin{lemma}\label{lemma_bplp}
\begin{enumerate}
\item[a)] Let $f$ be a b.p.l.p.~\Ref{eq_block_linearization}
satisfying~\Ref{eq_block_condition}. Then there is a tableau with substitution
$(T,(X_1,\ldots,X_s))$ such that ${\F}_T(X_1,\ldots,X_s)=\pm f$.

\item[b)] For every tableau with substitution $(T,(X_1,\ldots,X_s))$ of dimension
$\un{n}$ there is a b.p.l.p.~$f$ satisfying~\Ref{eq_block_condition} such that
${\F}_T(X_1,\ldots,X_s)=\pm f$.
\end{enumerate}
\end{lemma}
\begin{proof}
{\bf a)} A tableau with substitution $(T,(X_1,\ldots,X_s))$ of dimension $\un{n}$ is
constructed as follows. Its arrows are $a_{jr}$, where $a_{jr}$ goes from the $p_j$-th to
the $q_j$-th column for $1\leq j\leq s$ and $1\leq r\leq k_j$.
Condition~\Ref{eq_block_condition} guarantees that for any $1\leq i\leq m$ the total
number of arrows that begin or end in the $i$-th column is $n_i$. Complete the
construction by setting $\ovphi{a_{jr}}=j$. Note that $(T,(X_1,\ldots,X_s))$ is not
uniquely determined by $f$.


{\bf b)} Consider $a_1,\ldots,a_s\in T$ such that $\ovphi{a_1}=1,\ldots,\ovphi{a_s}=s$.
Then
$$\F_T(X_1,\ldots,X_s)=\pm\P_{k_1,\ldots,k_s}(X_1^{a_1'',a_1'},\ldots,X_s^{a_s'',a_s'}),$$
where $k_j=\#\{a\in T\,|\,\ovphi{a}=j\}$ for any $1\leq j\leq s$.
\end{proof}

\example\label{ex1} {\bf 1.} Let $n$ be even and $(T,(X_1,\ldots,X_s))$ be a tableau with
substitution of dimension $(n)$, where $T$ is
$$\begin{array}{|c|}
\hline a_1\vspace{-0.2cm}\\
\vspace{-0.3cm}\downarrow \\
\hline \\
\hline \vdots\\
\hline \! a_{n/2}\!\!\vspace{-0.2cm}\\
\vspace{-0.3cm}\downarrow \\
\hline \\
\hline
\end{array}\quad.
$$
In other words, arrows of $T$ are $a_1,\ldots,a_{n/2}$, where
$'a_i=2i$, $''a_i=2i-1$, $a_i'=a_i''=1$ for $1\leq i\leq n/2$, and
$X_1,\ldots,X_s$ are $n\times n$ matrices.

If $s=1$, then~\Ref{eq_P} implies $\F_T(X_1)=\P(X_1)$ for $K=\QQ$
and consequently for an arbitrary $K$.

If $s>1$, then by~\Ref{eq_lin_P} %
$${\F}_T(X_1,\ldots,X_s)={\P}_{k_1,\ldots,k_s}(X_1,\ldots,X_s),$$ %
where $k_j=\#\{a\in T\,|\,\ovphi{a}=j\}$ for any $1\leq j\leq s$,
is a partial linearization of the pfaffian.
\smallskip

{\bf 2.} For $n\times n$ matrices $X_1,\ldots,X_s$ let $(T,(X_1,\ldots,X_s))$ be the
tableau with substitution of dimension $(n,n)$, where $T$ is
$$\begin{array}{|c|c|}
\hline a_1-\!\!\!\! &\!\!\!\rightarrow\;\,\\
\hline \vdots &\vdots\\
\hline a_n-\!\!\!\! &\!\!\!\rightarrow\;\,\\
\hline
\end{array}\quad.
$$
If $s=1$, then the formula
\begin{eq}\label{eq_det}
\det(X)=\frac{1}{n!}\sum_{\pi_1,\pi_2\in
\Symmgr_n}\sign(\pi_1)\sign(\pi_2) \prod_{i=1}^n
x_{\pi_1(i),\pi_2(i)},
\end{eq}
which is valid over $\QQ$, implies the equality $\F_T(X_1)=\det(X_1)$ over every $K$. For
$s>1$ the expression $\F_T(X_1,\ldots,X_s)$ is a partial linearization of the
determinant.
\smallskip

{\bf 3.} If $(T,(X_1,X_2))$ is a tableau with substitution from part~2, where for $1\leq
k\leq n$ we have $\ovphi{a_1}=\cdots=\ovphi{a_k}=1$,
$\ovphi{a_{k+1}}=\cdots=\ovphi{a_n}=2$, $X_1=X$, and $X_2=E$ is the identity $n\times n$
matrix, then ${\F}_T(X,E)=\si_k(X)$.

{\bf 4.} Suppose that $t,s,r\in\NN$, and $X$, $Y$, $Z$,
respectively, are matrices of dimensions $(t+2s)\times(t+2r)$,
$(t+2s)\times(t+2s)$, $(t+2r)\times(t+2r)$, respectively. Let $T$
be a tableau of dimension $(t+2s,t+2r)$ that in the case $s=r$ is
depicted as
$$\begin{array}{|c|c|}
\hline a_1-\!\!\!\! &\!\!\!\rightarrow\;\,\\
\hline \vdots &\vdots\\
\hline a_t-\!\!\!\! &\!\!\!\rightarrow\;\,\\
\hline b_1\vspace{-0.2cm}&c_1\\
\vspace{-0.3cm}\downarrow& \downarrow\\
\hline &\\
\hline \vdots&\vdots\\
\hline b_{s}\vspace{-0.2cm}&c_{r}\\
\vspace{-0.3cm}\downarrow&\downarrow \\
\hline &\\
\hline
\end{array}\quad,
$$
otherwise define $T$ analogously. Let $(T,(X,Y,Z))$ be a tableau with substitution, where
$\ovphi{a_i}=1$, $\ovphi{b_j}=2$, and $\ovphi{c_k}=3$ for $1\leq i\leq t$, $1\leq j\leq
s$, and $1\leq k\leq r$. Then $\F_T(X,Y,Z)=\pm\DP_{s,r}(X,Y,Z)$, where $\DP_{s,r}$ was
introduced in Section~3 of~\cite{LZ1}.
\bigskip

The main result of the paper is the following decomposition
formula.

\begin{theo}\label{theo_short}
{\bf (Decomposition formula: short version)}\\
Let $(T,(X_1,\ldots,X_s))$ be a tableau with substitution of dimension $\un{n}\in\NN^m$.
Let $1\leq q_1<q_2\leq m$, $n_{q_1}=n_{q_2}$, and the vector $\un{d}\in\NN^{m-2}$ be
obtained from $\un{n}$ by eliminating the $q_1$-th and the $q_2$-th coordinates. Then
$\F_T(X_1,\ldots,X_s)$ is a polynomial in $\F_D(Y_1,\ldots,Y_l)$ and $\si_{k}(h)$, where
\begin{enumerate}
\item[$\bullet$] $(D,(Y_1,\ldots,Y_l))$ ranges over tableaux with substitutions of
dimension $\un{d}$ such that $Y_1,\ldots,Y_l$ are products of matrices
$X_1,\ldots,X_s,X_1^t,\ldots,X_s^t$;

\item[$\bullet$] $h$ ranges over products of matrices
$X_1,\ldots,X_s,X_1^t,\ldots,X_s^t$;

\item[$\bullet$] $k$ ranges over $[1,n_{q_1}]$. %
\end{enumerate} Moreover,
coefficients of this polynomial belong to the image of $\ZZ$ in
$K$ under the natural homomorphism.
\end{theo}

\section{Decomposition formula}\label{section_decomp_form}
Consider a tableau with substitution $(T,(X_1,\ldots,X_s))$ of dimension
$\un{n}\in\NN^m$. In what follows we assume that $n_{q_1}=n_{q_2}=n$ for $1\leq
q_1<q_2\leq m$. 



Denote by $\M$ the monoid freely generated by letters
$x_1,x_2,\ldots$, $x_1^t,x_2^t,\ldots$ Fix some lexicographical
order $<$ on $\M$. Let $\M_T$ be the submonoid of $\M$, generated
by $x_1,\ldots,x_s$, $x_1^t,\ldots,x_s^t$. For short, we will
write $1,\ldots,s,1^t,\ldots,s^t$ instead of
$x_1,\ldots,x_s,x_1^t,\ldots,x_s^t$. Given $a\in T$, we consider
$\ovphi{a}\in\{1,\ldots,s\}$ as an element of $\M_T$.

For $u\in\M_T$ define the matrix $X_u$ by the following rules:
\begin{enumerate}
\item[$\bullet$] $X_{j^t}=X_j^t$ for any $1\leq j\leq s$;

\item[$\bullet$] $X_{vw}=\left\{%
\begin{array}{cl}
X_v X_w,&\text{if the product of these matrices is well defined} \\
0,& \text{otherwise} \\
\end{array}
\right.$ %
\\for $v,w\in \M_T$.
\end{enumerate}

For an arrow $a\in T$ denote by $a^t$ the {\it transpose arrow},
i.e., by definition $(a^t)''=a'$, $(a^t)'=a''$, $''(a^t)={}'a$,
$'(a^t)={}''a$, $\ovphi{a^t}=\ovphi{a}^t\in \M_T$. Obviously,
$(a^t)^t=a$.

We write $a\stackrel{t}{\in} T$ if $a\in T$ or $a^t\in T$.

\bigskip 
\noindent{\bf Definition (of paths).} We say that
$a_1,a_2\stackrel{t}{\in}T$ are {\it successive} in $T$,  if
$a_1',a_2''\in\{q_1,q_2\}$, $a_1'\neq a_2''$, $'a_1={}''a_2$.

A word $a=a_1\cdots a_r$, where
$a_1,\ldots,a_r\stackrel{t}{\in}T$, is called a {\it path} in $T$
with respect to columns $q_1$ and $q_2$, if $a_i$, $a_{i+1}$ are
successive for any $1\leq i\leq r-1$. In this case by definition
$\ovphi{a}=\ovphi{a_1} \cdots \ovphi{a_r}\in \M_T$ and
$a^t=a_r^t\cdots a_1^t$ is a path in $T$; we denote
$a_r',{}'a_r,a_1'',{}''a_1$, respectively, by $a',{}'a,a'',{}''a$,
respectively. Since the columns $q_1$ and $q_2$ are fixed, we
usually refer to $a$ as a path in $T$.

A path $a_1\cdots a_r$ is {\it open} if both
$a_1'',a_r'\notin\{q_1,q_2\}$.

A path $a_1\cdots a_r$ is {\it closed} if $a_r,a_1$ are successive;   in particular,
$a_1'',a_r'\in\{q_1,q_2\}$. An element $a\stackrel{t}{\in} T$ with $a',a''\in\{q_1,q_2\}$
and $'a={}''a$ is also called a closed path.

A path is called {\it maximal} if it is either open or closed.

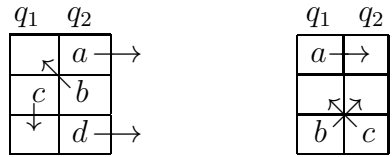
\example\label{ex2} Let %
$$\begin{array}{cc}
q_1& q_2\\
\end{array}
\qquad\qquad\qquad
\begin{array}{cc}
q_1&q_2\\
\end{array}
$$
\vspace{-.6cm}
$$\begin{array}{|c|c|} %
\hline & \vspace{-0.25cm}a\!\longrightarrow\hspace{-0.6cm}\\
\nwarrow\hspace{-0.6cm}
&\hspace{-0.5cm}\vspace{-0.25cm}\\
\hline \;c\vspace{-0.2cm}&\!b\\
\vspace{-0.3cm}\downarrow \\
\hline &d\!\longrightarrow\hspace{-0.6cm}\\
\hline
\end{array}
\qquad\qquad\qquad
\begin{array}{|c|c|} %
\hline a-\!\!\!\! &\!\!\!\!\rightarrow\\
\hline & \vspace{-0.25cm}\\
\nwarrow\hspace{-0.6cm}
&\hspace{-0.5cm}\vspace{-0.25cm}\nearrow\\
\hline b&\,c\\
\hline
\end{array}
$$%
be fragments of two tableaux, where in both cases we have drawn columns $q_1$, $q_2$. In
the first case we have an open path $a^tb^tcd$; and in the second case we have closed
paths: $a$, $bc^t$.
\bigskip

Consider words $a=a_1\cdots a_p$, $b=b_1\cdots b_q$, where
$a_1\cdots a_p,b_1\cdots b_q\stackrel{t}{\in}T$ (in particular,
both $a,b$ might be paths in $T$). We write
\begin{enumerate}
\item[$\bullet$] $a=b$ if $p=q$ and $a_1=b_1,\ldots,a_p=b_p$;

\item[$\bullet$] $a\stackrel{t}{=}b$ if $a=b$ or $a=b^t$;

\item[$\bullet$] $a\stackrel{c}{=}b$ if there is a cyclic
permutation $\pi\in \Symmgr_p$ such that $a_{\pi(1)}\cdots
a_{\pi(p)}=b$;

\item[$\bullet$] $a\stackrel{ct}{=}b$ if there is a cyclic
permutation $\pi\in \Symmgr_p$ such that $a_{\pi(1)}\cdots
a_{\pi(p)}\stackrel{t}{=}b$.
\end{enumerate}
Note that $\stackrel{t}{=}$, $\stackrel{c}{=}$, and $\stackrel{ct}{=}$ are equivalences.
We will use similar notations for elements of $\M_T$, for sets of paths, and for subsets
of $\M_T$. In the same fashion define the signs of inclusion $\stackrel{t}{\in}$,
$\stackrel{ct}{\in}$, and the sign of subtraction of sets $\stackrel{t}{\backslash}$. As
an example, if $s\geq 6$, then we have the following equivalence of subsets of $\M_T$
$$\{12^t,\, 3^t4,\, 56\}\,\stackrel{t}{\backslash}\,%
\{21^t\}\stackrel{t}{=}\{4^t3,\,56\}.$$ %
Note that
\begin{enumerate}
\item[$\bullet$] if $a$ is an open path and $a\stackrel{t}{=}b$,
then $b$ is also an open path;

\item[$\bullet$] if $a$ is a closed path and $a\stackrel{ct}{=}b$,
then $b$ is also a closed path;

\item[$\bullet$] if $a=a_1\cdots a_r$, where
$a_1,\ldots,a_r\stackrel{t}{\in}T$, is an open path, then
$\ovphi{a}\in\M_T$ is a primitive element (see
Section~\ref{section_prelim} for the definition), since otherwise
there is an $1<l\leq r$ such that $\ovphi{a_1}=\ovphi{a_l}$; hence
$a_1''={}a_l''$; on the other hand, $a_1''\notin\{q_1,q_2\}$ and
$a_l''\in\{q_1,q_2\}$; a contradiction.
\end{enumerate}

Denote by $T_{o}$ the set $\{b\,|\,b\text{ is an open path in
}T,\,(b'',{}''b)<(b',{}'b)\}$, where $<$ is the lexicographical
order on $\NN\times\NN$ introduced in
Section~\ref{section_prelim}. Obviously, $T_o$ is a set of
representatives of open paths in $T$ with respect to the
$\stackrel{t}{=}$-equivalence, i.e., $a\stackrel{t}{\in}T_{o}$ for
every open path $a$ in $T$, and $a\stackrel{t}{\neq}b$ for all
$a,b\in T_{o}$ with $a\neq b$. Let $T_{cl}$ be a set of
representatives of closed paths in $T$ with respect to the
$\stackrel{ct}{=}$-equivalence. Define %
$$\ovphi{T_{o}}=\Lus \ovphi{a}\,|\, a\in T_{o}\Rus,\;\; %
\ovphi{T_{cl}}=\Lus \ovphi{a} \,|\, a\in T_{cl}\Rus$$ %
(see Section~\ref{section_prelim} for the definition of a multiset). Then the expression
$$\tr(T_{cl})=\prod_{a\in T_{cl}}\tr(X_{\ovphi{a}}).$$ %
is well defined.

\bigskip 
\noindent {\bf Definition (of $T^{\tau}$ and $\widetilde{T}$).} %
{\bf a)} Let $\tau\in \Symmgr_{n}$. Permute the cells of the
$q_2$-th column of $T$ by $\tau$ and denote the resulting tableau
by $T^{\tau}$. The arrows of $T^{\tau}$ are $\{a^{\tau}\,|\,a\in
T\}$, where $\ovphi{a^{\tau}}=\ovphi{a}$, $(a^{\tau})''=a''$,
$(a^{\tau})'=a'$, and %
$$''(a^{\tau})=\left\{
\begin{array}{cl}
\tau(''a),& \text{if } a''=q_2\\
''a,&\text{otherwise}\\
\end{array}
\right. ,\quad%
'(a^{\tau})=\left\{
\begin{array}{cl}
\tau('a),&\text{if } a'=q_2\\
'a,&\text{otherwise}\\
\end{array}
\right..
$$
Obviously, $(T^{\tau},(X_1,\ldots,X_s))$ is a tableau with substitution.

{\bf b)} Define the tableau with substitution $(\widetilde{T},X_{\widetilde{T}})$ of
dimension $(n_1,\ldots,n_{q_1-1}$, $n_{q_1+1},\ldots,n_{q_2-1}$, $n_{q_2+1},\ldots,n_m)$,
where $X_{\widetilde{T}}=(Y_1,\ldots,Y_t)$ for some matrices $Y_1,\ldots,Y_t$ for
suitable $t>0$ as follows. If $q_1=m-1$, $q_2=m$, then $\widetilde{T}$ contains two types
of arrows:
\begin{enumerate}
\item[$\bullet$] $\widetilde{a}$, where $a\in T$, and both
$a',a''\notin\{q_1,q_2\}$;

\item[$\bullet$] $\widetilde{a}$, where $a$ is an open path from $T_{o}$. %
\end{enumerate}
In both cases the tail and the head of $\widetilde{a}$ coincides with the tail and the
head of $a$, respectively. Define $t>0$ and $\ovphi{u}$ for $u\in\widetilde{T}$ in such a
way that $\{\ovphi{u}\,|\,u\in \widetilde{T}\}=[1,t]$ and for all
$\widetilde{b},\widetilde{c}\in\widetilde{T}$ we have
\begin{itemize}
\item the equality $\ovphi{\widetilde{b}}=\ovphi{\widetilde{c}}$
holds if and only if $\ovphi{b}=\ovphi{c}$,

\item  the inequality
$\ovphi{\widetilde{b}}<\ovphi{\widetilde{c}}$ holds if and only if
$\ovphi{b}<\ovphi{c}$,
\end{itemize}
where $<$ is the lexicographical order on $\M$ fixed at the beginning of this section. We
set $Y_{\ovphi{\widetilde{a}}}=X_{\ovphi{a}}$. Then the tableau with substitution
$(\widetilde{T},X_{\widetilde{T}})$ is well defined.

For arbitrary $q_1$, $q_2$ the definition of
$(\widetilde{T},X_{\widetilde{T}})$ is obtained analogously using
an appropriate shift to define numbers
$\widetilde{a}',\widetilde{a}''$.

For short, we will omit parentheses and use the following
conventions%
$$\begin{array}{ccc}
(T^{\tau})_{o}=T^{\tau}_{o},& (T^{\tau})_{cl}=T^{\tau}_{cl},&\\
\widetilde{(T^{\tau})}=\widetilde{T}^{\tau},&
\widetilde{(a^{\tau})}=\widetilde{a}^{\tau},&\\
(\widetilde{T}^{\tau})_{o}=\widetilde{T}^{\tau}_{o},&
(\widetilde{T}^{\tau})_{cl}=\widetilde{T}^{\tau}_{cl},&\\
\end{array}
$$ %
for any $\tau\in \Symmgr_n$ and $a\in T$.

\example\label{ex3} {\bf 1.}
Let $T$ be the  tableau  %
$\begin{array}{|c|c|c|} %
\hline a-\!\!\!\!  &  \!\!\!\rightarrow  &  b\vspace{-0.2cm}\\
\vspace{-0.3cm}&& \downarrow\\
\hline c-\!\!\!\!  & \!\!\!\rightarrow  & \\
\hline
\end{array}\;$ %
of dimension $(2,2,2)$. Suppose $q_1=2$ and $q_2=3$. Then
$\widetilde{T}$ is %
$\begin{array}{|c|} %
\hline d\vspace{-0.2cm}\\
\vspace{-0.3cm}\downarrow \\
\hline \\
\hline
\end{array}\;$, %
where $d=\widetilde{abc^t}$.

\smallskip
{\bf 2.} Let $T$ be the tableau
$$\begin{array}{|c|c|c|}
\hline \vspace{-0.2cm} a & b-\!\!\!\! &\!\!\!\!\rightarrow\,\\
\vspace{-0.3cm}\downarrow && \\ 
\hline \vspace{-0.25cm}& &\\
\vspace{-0.25cm}\nearrow\!\!\!\!\!\!\!\!\!&\nearrow\!\!\!\!\!\!\!\!\! \\    
\hline c& d\\
\hhline{--~}
\end{array}
$$ %
of dimension $(3,3,2)$. Suppose $q_1=1$, $q_2=2$, and $\tau\in
\Symmgr_3$ is the cycle $(1,2,3)$. Then $T^{\tau}$ is
$$\begin{array}{|c|c|c|}
\hline \vspace{-0.2cm} a^{\tau} & d^{\tau} &\;\;\;\\
\vspace{-0.3cm}\downarrow &\searrow\hspace{-0.75cm}
&\hspace{-0.7cm}\nearrow \\ 
\hline \vspace{-0.25cm}& b^{\tau}&\\
\vspace{-0.25cm}& \\    
\hline c^{\tau}-\!\!\!\! &\!\!\!\!\rightarrow\\
\hhline{--~}
\end{array}\;,
$$
$\widetilde{T}$ is %
$\begin{array}{|c|} %
\hline e\vspace{-0.2cm}\\
\vspace{-0.3cm}\downarrow \\
\hline \\
\hline
\end{array}\;$, %
and $\widetilde{T}^{\tau}$ is %
$\begin{array}{|c|} %
\hline h\vspace{-0.2cm}\\
\vspace{-0.3cm}\downarrow \\
\hline \\
\hline
\end{array}\;$, %
where $e=\widetilde{b^tac^td}$, $h=\widetilde{(b^{\tau})^t(a^{\tau})^t d^{\tau}}$. In
particular, $\widetilde{T}$ and $\widetilde{T}^{\tau}$ are different.

\bigskip
\noindent{\bf Definition (of admissibility).} Let
$\un{\be}\in\NN^p$, $\un{\ga}\in\NN^q$, $\un{b}=(b_1,\ldots,b_p)$,
$\un{c}=(c_1,\ldots,c_q)$, where $b_i$, $c_j$  are primitive
elements of $\M_T$. Moreover, let $b_1,\ldots,b_p$ be pairwise
different with respect to $\stackrel{t}{=}$, i.e.,
$b_i\stackrel{t}{\neq}b_j$ for $i\neq j$; and let $c_1,\ldots,c_q$
be pairwise different with respect to $\stackrel{ct}{=}$. Then
$(\un{\be},\un{\ga},\un{b},\un{c})$ is called a {\it
$T$-quadruple}.

A $T$-quadruple $(\un{\be},\un{\ga},\un{b},\un{c})$ is called {\it
$T$-admissible} if for some
$\xi=\xi_{\un{\be},\un{\ga},\un{b},\un{c}}\in \Symmgr_n$ the
following equivalences of multisets hold: %
$$\ovphi{T^{\xi}_{o}}\stackrel{t}{=}%
\Lus\underbrace{b_1,\ldots,b_1}_{\be_1},\ldots,\underbrace{b_p,\ldots,b_p}_{\be_p}\Rus,\;\;%
\ovphi{T^{\xi}_{cl}}\stackrel{ct}{=}%
\Lus\underbrace{c_1,\ldots,c_1}_{\ga_1},\ldots,\underbrace{c_q,\ldots,c_q}_{\ga_q}\Rus.$$
We write
$(\un{\be^0},\un{\ga^0},\un{b^0},\un{c^0})\stackrel{ct}{=}(\un{\be},\un{\ga},\un{b},\un{c})$
and say that these quadruples are equivalent if and only if %
$$
\Lus\underbrace{b_1^0,\ldots,b_1^0}_{\be_1^0},\ldots,
\underbrace{b_p^0,\ldots,b_p^0}_{\be_p^0}\Rus%
\stackrel{t}{=}%
\Lus\underbrace{b_1,\ldots,b_1}_{\be_1},\ldots,\underbrace{b_p,\ldots,b_p}_{\be_p}\Rus,$$
$$ %
\Lus\underbrace{c_1^0,\ldots,c_1^0}_{\ga_1^0},\ldots,
\underbrace{c_q^0,\ldots,c_q^0}_{\ga_q^0}\Rus%
\stackrel{ct}{=}%
\Lus\underbrace{c_1,\ldots,c_1}_{\ga_1},\ldots,\underbrace{c_q,\ldots,c_q}_{\ga_q}\Rus.$$
If $(\un{\be^0},\un{\ga^0},\un{b^0},\un{c^0})\stackrel{ct}{=}
(\un{\be},\un{\ga},\un{b},\un{c})$ and
$(\un{\be},\un{\ga},\un{b},\un{c})$ is $T$-admissible, then the
quadruple $(\un{\be^0},\un{\ga^0},\un{b^0},\un{c^0})$ also has the
same property, since we can take
$\xi_{\un{\be^0},\un{\ga^0},\un{b^0},\un{c^0}}=\xi_{\un{\be},\un{\ga},\un{b},\un{c}}$.
Denote by $Q_T$ a set of representatives of $T$-admissible
quadruples with respect to the $\stackrel{ct}{=}$-equivalence.

\begin{theo}{\bf(Decomposition formula).}\label{theo_decomp} 
Let  $1\leq q_1<q_2\leq m$ and $n_{q_1}=n_{q_2}=n$. Then for a tableau with substitution
$(T,(X_1,\ldots,X_s))$ of dimension $\un{n}$ we have
$${\F}_T(X_1,\ldots,X_s)=\sum_{(\un{\be},\un{\ga},\un{b},\un{c})\in Q_T}  %
\sign(\xi_{\un{\be},\un{\ga},\un{b},\un{c}}) %
{\F}_{\widetilde{T}^{\xi_{\un{\be},\un{\ga},\un{b},\un{c}}}} %
(X_{\widetilde{T}^{\xi_{\un{\be},\un{\ga},\un{b},\un{c}}}})\prod_{j=1}^{\#\un{\ga}}
\si_{\ga_i}(X_{c_j}).$$ %
\end{theo}
\bigskip

Neither the permutation $\xi_{\un{\be},\un{\ga},\un{b},\un{c}}$ nor its sign are unique
for a representative of the $\stackrel{ct}{=}$-equivalence class of a $T$-admissible
quadruple $(\un{\be},\un{\ga},\un{b},\un{c})$. However, it follows from part~b) of
Lemma~\ref{lemma_coef}, which is formulated below, that
$\sign(\xi_{\un{\be},\un{\ga},\un{b},\un{c}})
\F_{\widetilde{T}^{\xi_{\un{\be},\un{\ga},\un{b},\un{c}}}}
(X_{\widetilde{T}^{\xi_{\un{\be},\un{\ga},\un{b},\un{c}}}})$ does not depend on a
representative of the $\stackrel{ct}{=}$-equivalence class. Thus the right hand side of
the decomposition formula is well defined.

Sections~\ref{section_decomp_F0},~\ref{section_coef},
and~\ref{section_tr} are devoted to a proof of this theorem. Since
the right hand side of the decomposition formula is a polynomial
with integer coefficients in $\ZZ$, if the characteristic of $K$
is zero, or in $\ZZ/p\ZZ$, if the characteristic of $K$ is $p>0$;
it is enough to prove the theorem in the case $K=\QQ$. To simplify
notation we assume that $q_1=m-1$, $q_2=m$.

The proof is organized as follows. In Section~\ref{section_decomp_F0} we prove
Lemma~\ref{lemma_decomp} for $\F_T^0(X_1,\ldots,X_s)$, which is a statement analogous to
Theorem~\ref{theo_decomp}. The short version of the decomposition formula (see
Theorem~\ref{theo_short}) follows from Lemma~\ref{lemma_decomp} in the case of
characteristic zero. In order to prove Theorem~\ref{theo_short} as well as
Theorem~\ref{theo_decomp} over the field of arbitrary characteristic we, working over
$\QQ$, should divide the formula in Lemma~\ref{lemma_decomp} by $c_T$. This problem is
solved by rewriting $\F_T(X_1,\ldots,X_s)$ in a more suitable form (see
formula~\Ref{eq_coef} in Section~\ref{section_coef}) and by applying
Lemma~\ref{lemma_zeta} from Section~\ref{section_tr}.

\section{Decomposition of ${\F}_T^0(X_1,\ldots,X_s)$}\label{section_decomp_F0}
\begin{lemma}\label{lemma_decomp}
Let $\un{n}\in\NN^m$, $1\leq q_1<q_2\leq m$, $n_{q_1}=n_{q_2}=n$, and
$(T,(X_1,\ldots,X_s))$ be a tableau with substitution of dimension $\un{n}$. Then
$$\F_T^0(X_1,\ldots,X_s)=\sum_{\tau\in \Symmgr_n} \sign(\tau)
\F^0_{\widetilde{T}^{\tau}}(X_{\widetilde{T}^{\tau}})\tr(T_{cl}^{\tau}).$$
\end{lemma}
\begin{proof}
Without loss of generality we can assume that $q_1=m-1$, $q_2=m$. Define %
$$A=\{ a\in T\,|\,a'>m-2\text{ or }a''>m-2\}.$$ %
For fixed $\pi_1\in \Symmgr_{n_1},\ldots,\pi_{m-2}\in
\Symmgr_{n_{m-2}}$ define %
$$f_0=\sum_{\pi_{m-1},\pi_m\in \Symmgr_n}\sign(\pi_{m-1})\sign(\pi_m)\prod_{a\in A}
(X_{\ovphi{a}})_{\pi_{a''}(''a),\pi_{a'}('a)},$$
$$f_{T}=\prod_{b\in T_o}
(X_{\ovphi{b}})_{\pi_{b''}(''b),\pi_{b'}('b)}\,
\tr(T_{cl}).$$ %
We will show that
\begin{eq}\label{eq_f0}
f_0=\sum_{\tau\in \Symmgr_n} \sign(\tau) f_{T^{\tau}}.
\end{eq}

For an open path $b=b_1\cdots b_r\in T_o$ with
$b_1,\ldots,b_r\stackrel{t}{\in}T$ numbers
$'b_1,\ldots,'b_{r-1}$ are pairwise different. It is not difficult to see that %
$$(X_{\ovphi{b}})_{\pi_{b''}(''b),\pi_{b'}('b)}=
\sum_{1\leq h('b_1),\ldots,h('b_{r-1})\leq n} B(h('b_1),\ldots,h('b_{r-1})), $$%
where $B(h('b_1),\ldots,h('b_{r-1}))$ stands for
$$(X_{\ovphi{b_1}})_{\pi_{b_1''}(''b_1),h('b_1)}\,%
(X_{\ovphi{b_2}})_{h(''b_2),h('b_2)}\cdots %
(X_{\ovphi{b_r}})_{h(''b_r),\pi_{b_r'}('b_r)}.$$ %
Similarly, for a closed path $c=c_1\cdots c_l\in T_{cl}$ with
$c_1,\ldots,c_l\stackrel{cl}{\in}T$ numbers
$'c_1,\ldots,'c_{l}$ are pairwise different. Therefore %
$$\tr(X_{\ovphi{c}})=\sum_{1\leq h('c_1),\ldots,h('c_{l})\leq n}
C(h('c_1),\ldots,h('c_l)),$$
where $C(h('c_1),\ldots,h('c_l))$ stands for %
$$(X_{\ovphi{c_1}})_{h(''c_1),h('c_1)}\, %
(X_{\ovphi{c_2}})_{h(''c_2),h('c_2)}\cdots
(X_{\ovphi{c_l}})_{h(''c_{l}),h('c_l)}.$$ %
Summarizing, we obtain %
$$f_T=\prod_{b_1\cdots b_r\in T_o}\sum_{1\leq h('b_1),\ldots,h('b_{r-1})\leq n}
B(h('b_1),\ldots,h('b_{r-1})) $$ %
$$*\prod_{c_1\cdots c_l\in
T_{cl}}\sum_{1\leq h('c_1),\ldots,h('c_l)\leq n}
C(h('c_1),\ldots,h('c_l)).$$ %
Since elements of the set $\{'b_i,{}'c_j\,|\,r,l>0,\,b_1\cdots
b_r\in T_o,\, c_1\cdots c_l\in T_{cl},\,1\leq i\leq r-1,\,1\leq
j\leq l\}$ are pairwise different
and their union equals $[1,n]$, we infer %
$$f_T=\sum_{h:[1,n]\to[1,n]} \;\prod_{b_1\cdots b_r\in T_o,\,c_1\cdots c_l\in T_{cl}} %
B(h('b_1)\cdots h('b_{r-1}))\, C(h('c_1)\cdots h('c_l)).$$ %
For an $a\stackrel{t}{\in}T$ and a function $h:[1,n]\to[1,n]$
we set: %
$$(h,{}a)^1=\left\{
\begin{array}{cl}
h('a),& \text{if } a'\in\{q_1,q_2\} \\
\pi_{a'}('a),& \text{otherwise}
\end{array} \right.,$$ %
$$(h,{}a)^2=\left\{
\begin{array}{cl}
h(''a),& \text{if } a''\in\{q_1,q_2\} \\
\pi_{a''}(''a),& \text{otherwise}
\end{array} \right..$$ %
Recall that $\pi_1\in \Symmgr_{n_1},\ldots,\pi_{m-2}\in
\Symmgr_{n_{m-2}}$ have been fixed. Since
$(X_{\ovphi{a^t}})_{(h,{a^t})^2,(h,{a^t})^1}=
(X_{\ovphi{a}})_{(h,{a})^2,(h,{a})^1}$ for all $a\in T$, we can rewrite $f_T$ in a form %
$$f_T=\sum_{h:[1,n]\to[1,n]}\;\prod_{a\in A}
(X_{\ovphi{a}})_{(h,{a})^2,(h,{a})^1}.$$ %
Denote by $H$ the set of functions $h:[1,n]\to[1,n]$ that are not
bijections. Since bijections $h:[1,n]\to[1,n]$ are
in one to one correspondence with permutations from $\Symmgr_n$, we obtain that %
$$ \sum_{\tau\in \Symmgr_n} \sign(\tau)f_{T^{\tau}}$$ %
$$=\sum_{\tau,\pi\in \Symmgr_n} \sign(\tau)\prod_{a\in A}
(X_{\ovphi{a}})_{(\pi,{}a^{\tau})^2,(\pi,{}a^{\tau})^1}+ %
\sum_{\tau\in \Symmgr_n, h\in H} \sign(\tau)\prod_{a\in A}
(X_{\ovphi{a}})_{(h,{}a^{\tau})^2,(h,{}a^{\tau})^1}.
$$ %
Denote the first and the second summands of the last expression by
$f_1$ and $f_2$, respectively.

Substituting $\pi_{m-1}^{-1}\pi_{m}$ for $\tau$ and $\pi_{m-1}$
for $\pi$ in
$f_1$ we get %
$$f_1=\sum_{\pi_{m-1},\pi_{m}\in \Symmgr_n} \sign(\pi_{m-1})\sign(\pi_{m})\prod_{a\in A}
(X_{\ovphi{a}})_{(\pi_{m-1},{}a^{\pi_{m-1}^{-1}\pi_{m}})^2,
(\pi_{m-1},{}a^{\pi_{m-1}^{-1}\pi_{m}})^1}.
$$ %
The following calculations show that
$(\pi_{m-1},{}a^{\pi_{m-1}^{-1}\pi_{m}})^1=\pi_{a'}({}'a)$.
\begin{enumerate}
\item[$\bullet$] If $a'=q_1=m-1$, then
$(\pi_{m-1},{}a^{\pi_{m-1}^{-1}\pi_{m}})^1=
\pi_{m-1}('a^{\pi_{m-1}^{-1}\pi_{m}})=\pi_{m-1}('a)$.

\item[$\bullet$] If $a'=q_2=m$, then
$(\pi_{m-1},{}a^{\pi_{m-1}^{-1}\pi_{m}})^1=
\pi_{m-1}('a^{\pi_{m-1}^{-1}\pi_{m}})=\pi_{m-1}\circ
\pi_{m-1}^{-1}\circ\pi_{m}('a)=\pi_{m}('a)$.

\item[$\bullet$] If $a'<m-1$, then
$(\pi_{m-1},{}a^{\pi_{m-1}^{-1}\pi_{m}})^1=
\pi_{k}('a^{\pi_{m-1}^{-1}\pi_{m}})=\pi_{a'}('a)$, where
$k=(a^{\pi_{m-1}^{-1}\pi_{m}})'=a'$.
\end{enumerate} %
In the same way we obtain
$(\pi_{m-1},{}a^{\pi_{m-1}^{-1}\pi_{m}})^2=\pi_{a''}({}''a)$.
Therefore, $f_1=f_0$.

The set $H$ can be represented as the disjoint union
$H=\bigsqcup\limits_{1\leq p<q\leq n} H_{pq}$ of sets %
$$H_{pq}=\{h\in H\,|\,h(p)=h(q)
{\rm\;and\;} h(i)\neq h(j) {\rm\;for\;every\;}(i,j)<(p,q)\},$$%
where $<$ is the lexicographical order on $\NN\times\NN$ defined
in Section~\ref{section_prelim}. We have
$$f_2=\sum_{1\leq p<q\leq n}\; \sum_{h\in H_{pq}}\left( %
 \sum_{\tau\in \Symmgr_n,\, \tau^{-1}(p)<\tau^{-1}(q)}\sign(\tau)\prod_{a\in
 A}
 (X_{\ovphi{a}})_{(h,{a^{\tau}})^2,(h,{}a^{\tau})^1}- \right.$$ %
 $$ \left. \sum_{\tau\in \Symmgr_n, \tau^{-1}(p)<\tau^{-1}(q)} \sign(\tau)
 \prod_{a\in A}
 (X_{\ovphi{a}})_{(h,{a^{(p,q)\circ\tau}})^2,(h,{}a^{(p,q)\circ\tau})^1}
\right),$$%
where, as usual, $(p,q)$ stands for a transposition from
$\Symmgr_n$. Compute $(h,{}a^{(p,q)\circ\tau})^1$ for $h\in
H_{pq}$.
\begin{enumerate}
\item[$\bullet$] If $a'=m-1$, then
$(h,{}a^{(p,q)\circ\tau})^1=h('a)=h('a^{\tau})$.

\item[$\bullet$] If $a'=m$, then
$(h,{}a^{(p,q)\circ\tau})^1=h\circ
(p,q)\circ\tau('a)=h\circ\tau('a)=h('a^{\tau})$.

\item[$\bullet$] If $a'<m-1$, then
$(h,{}a^{(p,q)\circ\tau})^1=\pi_{(a^{(p,q)\circ\tau})'}('a^{(p,q)\circ\tau})=\pi_{a'}('a)=
\pi_{a'}('a^{\tau})$. %
\end{enumerate}
Hence $(h,{}a^{(p,q)\circ\tau})^1=(h,{}a^{\tau})^1$ and,
similarly, $(h,{}a^{(p,q)\circ\tau})^2=(h,{}a^{\tau})^2$. This
implies $f_2=0$ and proves~\Ref{eq_f0}.

Rewrite $\F^0_T(X_1,\ldots,X_s)$ in a form %
$$\sum_{\pi_1\in \Symmgr_{n_1},\ldots,\pi_{m-2}\in \Symmgr_{n_{m-2}}}
\!\sign(\pi_1)\cdots\sign(\pi_{m-2})\! \prod_{a\in T, a\not\in A}
(X_{\ovphi{a}})_{\pi_{a''}(''a),\pi_{a'}('a)} \cdot f_0.
$$%
Applying~\Ref{eq_f0}, and taking into account that %
$$\prod_{b\in T_o^{\tau}}
(X_{\ovphi{b}})_{\pi_{b''}(''b),\pi_{b'}('b)}=%
\prod_{a=\widetilde{b},\,b\in T_{o}^{\tau}}%
(X_{\ovphi{a}})_{\pi_{a''}(''a),\pi_{a'}('a)},$$ %
for all $\tau\in\Symmgr_n$ we complete the proof.
\end{proof}

\section{Coefficients}\label{section_coef}
Continuing the proof of Theorem~\ref{theo_decomp} we keep the
assumptions and notations from previous sections. The aim of this
section is to prove formula~\Ref{eq_coef} (see below).

Let $(\un{\be},\un{\ga},\un{b},\un{c})$ be a $T$-admissible
quadruple with $\#\un{b}=p$ and $\#\un{c}=q$. Denote by
$\Symmgr_{\un{\be},\un{\ga},\un{b},\un{c}}$ the subset of
$\Symmgr_n$ containing $\tau\in \Symmgr_n$ if and only if %
$$\begin{array}{ccc}
\ovphi{T^{\tau}_{o}}&\stackrel{t}{=}&%
\Lus\underbrace{b_1,\ldots,b_1}_{\be_1},\ldots,\underbrace{b_p,\ldots,b_p}_{\be_p}\Rus
\;\text{ and}\\
\ovphi{T_{cl}^{\tau}} & \stackrel{ct}{=}& %
\Lus c_1^{\la_{11}},\ldots,c_1^{\la_{1l_1}},\ldots,
c_q^{\la_{q1}},\ldots, c_q^{\la_{ql_q}}\Rus,\\
\end{array}$$ %
for some positive integers $\la_{i1},\ldots,\la_{il_i}$ satisfying
$\la_{i1}+\cdots+\la_{il_i}=\ga_i$ $(1\leq i\leq q)$, where $c_i^{\la_{ij}}$ stands for
$c_i\cdots c_i$ ($\la_{ij}$ times). Notice that $\ovphi{T_o^{\tau}}$ and
$\ovphi{T_{cl}^{\tau}}$ are multisets (see Section~\ref{section_decomp_form}). For
$\tau\in \Symmgr_n$ we set
$$\Symmgr_{\tau}=\{\pi\in \Symmgr_n\,|\,\ovphi{T_o^{\pi}}\stackrel{t}{=}\ovphi{T_o^{\tau}},\,%
\ovphi{T_{cl}^{\pi}}\stackrel{ct}{=}\ovphi{T_{cl}^{\tau}}\}.$$ %
The following properties are immediate consequences of the
definitions.

\begin{lemma}\label{lemma_abcd}
\begin{enumerate}
\item[a)] For all $\pi,\tau\in \Symmgr_n$ we have
$\Symmgr_{\pi}=\Symmgr_{\tau}$ or $\Symmgr_{\pi}\cap
\Symmgr_{\tau}=\emptyset$;

\item[b)] $\tau\in \Symmgr_{\tau}$ for all $\tau\in \Symmgr_n$;

\item[c)] If $\tau\in \Symmgr_{\un{\be},\un{\ga},\un{b},\un{c}}$,
then $\Symmgr_{\tau}\subset
\Symmgr_{\un{\be},\un{\ga},\un{b},\un{c}}$;

\item[d)]
$\Symmgr_n=\bigsqcup\limits_{(\un{\be},\un{\ga},\un{b},\un{c})\in
Q_T} \Symmgr_{\un{\be},\un{\ga},\un{b},\un{c}}$.
\end{enumerate} %
\end{lemma}
\bigskip

\noindent Denote by $\ov{\Symmgr}_{\un{\be},\un{\ga},\un{b},\un{c}}$ a
subset of $\Symmgr_{\un{\be},\un{\ga},\un{b},\un{c}}$ such that
\begin{eq}\label{eq_ovS}
\Symmgr_{\un{\be},\un{\ga},\un{b},\un{c}}=\bigsqcup_{\tau\in
\ov{\Symmgr}_{\un{\be},\un{\ga},\un{b},\un{c}}}\Symmgr_{\tau}. %
\end{eq}
The existence of $\ov{\Symmgr}_{\un{\be},\un{\ga},\un{b},\un{c}}$ follows from
parts~a),~b),~c) of Lemma~\ref{lemma_abcd}.

Rewrite the coefficient $c_T=\prod_{j=1}^{s}\#\{a\in
T\,|\,\ovphi{a}=j\}!$ as
$c_T=c_T^{(q_1,q_2)}c_T^{(rest)}$, where%
$$c^{(q_1,q_2)}_T= %
\prod_{j=1}^s\#\{a\in T\,|\,\ovphi{a}=j,\, a'>m-2\;{\rm or}\;a''>m-2 \}!%
{\;\;\rm and}$$ %
$$c^{(rest)}_T= %
\prod_{j=1}^s\#\{a\in T\,|\,\ovphi{a}=j,\, a',a''\leq m-2\}!.$$ %
Given $\tau\in \Symmgr_n$, we have
$c_{\widetilde{T}^{\tau}}=c^{(rest)}_{T}c^{(\tau)}_T$, where %
$$c^{(\tau)}_T=%
\prod_{j=1}^t\#\{b\in T_{o}^{\tau}\,|\,\ovphi{\widetilde{b}}=j\}!$$ %
for $X_{\widetilde{T}}=(Y_1,\ldots,Y_t)$. According to part~d) of
Lemma~\ref{lemma_abcd}, for every $\tau\in \Symmgr_n$ there is a
unique (up to the $\stackrel{ct}{=}$-equivalence) quadruple
$(\un{\be},\un{\ga},\un{b},\un{c})\in Q_T$ such that $\tau\in
\Symmgr_{\un{\be},\un{\ga},\un{b},\un{c}}$. Therefore %
$$\ovphi{T_{cl}^{\tau}}\stackrel{ct}{=}
\Lus c_1^{\la_{11}},\ldots,c_1^{\la_{1l_1}},\ldots,
c_q^{\la_{q1}},\ldots, c_q^{\la_{ql_q}}\Rus$$ %
for some partitions $\un{\la}_i=(\la_{i1},\ldots,\la_{il_i})\vdash
\ga_i$, where $1\leq i\leq q$. Denote the corresponding multi-partition by  %
\begin{eq}\label{eq_la}
\un{\la}(\tau)=(\un{\la_1},\ldots,\un{\la_q})\vdash \un{\ga}.
\end{eq}

For a vector $\un{\la}=(\la_1,\ldots,\la_l)\in\NN^l$ define
$\un{\la}!=\la_1!\cdots\la_l!$
and %
\begin{eq}\label{eq_c_la}
c(\un{\la})=\la_1\cdots\la_l\cdot\prod_{k>0}\#\{1\leq i\leq
l\,|\,\la_i=k \}!.
\end{eq} %
For a multivector $\un{\la}=(\un{\la}_1,\ldots,\un{\la}_q)$ with
$\un{\la}_i\in\NN^{l_i}$ we set
$c(\un{\la})=c(\un{\la}_1)\cdots(\un{\la}_q)$.

\smallskip \noindent{\bf Definition (of the $\sim$-equivalence).}
For $a_1,\ldots,a_p,b_1,\ldots,b_q\stackrel{t}{\in}T$ and $\tau\in \Symmgr_n$ consider
words $a=a_1\cdots a_p$, $b=b_1\cdots b_q$, and $b^{\tau}=b_1^{\tau}\cdots b_q^{\tau}$.
If $a=b$, i.e., $p=q$ and $a_1=b_1,\ldots,a_p=b_p$, then we write $a\sim b^{\tau}$ and
call $a$, $b^{\tau}$ {\it equivalent}. Analogously, if $a\stackrel{t}{=}b$ or
$a\stackrel{ct}{=}b$, respectively, then we write $a\stackrel{t}{\sim}b^{\tau}$ or
$a\stackrel{ct}{\sim}b^{\tau}$, respectively. We will also use similar notation for
multisets.

\example\label{ex35} Let $T$ be the tableau  %
$$
\begin{array}{|c|c|c|}
\hline \vspace{-0.2cm} a & b-\!\!\!\! &\!\!\!\!\rightarrow\,\\
\vspace{-0.3cm}\downarrow && \\ 
\hline \vspace{-0.25cm}& &\\
\vspace{-0.25cm}\nearrow\!\!\!\!\!\!\!\!\!&\nearrow\!\!\!\!\!\!\!\!\! \\    
\hline \vspace{-0.2cm} c& d&e\\
\vspace{-0.3cm} && \downarrow\\ 
\hline f-\!\!\!\! &\!\!\!\!\rightarrow\,&\\
\hline
\end{array}
$$
of dimension $(4,4,4)$. Suppose $q_1=1$, $q_2=2$, and $\tau\in \Symmgr_4$ is such a
permutation that $T_{o}^{\tau}=\Lus (d^{\tau})^t (a^{\tau})^t b^{\tau}\Rus$ and
$T_{cl}^{\tau}=\Lus c^{\tau}f^{\tau}\Rus$. Then $T_{o}^{\tau} \sim \Lus d^t a^t b\Rus$,
$T_{cl}^{\tau}\sim\Lus c f\Rus$, and $\tau$ is uniquely determined by these data, namely,
$\tau$ is the cycle $(2,4,3)$ (see Lemma~\ref{lemma_remark} and its
proof). In particular, $T^{\tau}$ is %
$$\begin{array}{|c|c|c|} %
\hline a^{\tau}\vspace{-0.2cm}  &   b^{\tau}-\!\!\!\!  & \!\!\!\rightarrow\\
\downarrow&\vspace{-0.3cm}& \\
\hline &  d^{\tau}-\!\!\!\!  &\!\!\!\rightarrow  \\
\hline \vspace{-0.2cm}  c^{\tau} &\;\;&e^{\tau} \\
\vspace{-0.3cm}\;\searrow\hspace{-0.75cm}
&\hspace{-0.75cm}\nearrow&\downarrow  \\ 
\hline  f^{\tau}&&\\
\hline
\end{array} %
$$

\begin{lemma}\label{lemma_remark}
Suppose $\pi_1,\pi_2\in \Symmgr_n$. Then $\pi_1=\pi_2$ if and only if
$T_{o}^{\pi_1}\stackrel{t}{\sim}T_{o}^{\pi_2}$ and
$T_{cl}^{\pi_1}\stackrel{ct}{\sim}T_{cl}^{\pi_2}$.
\end{lemma}
\begin{proof}\textbf{1.} If $\pi_1=\pi_2$, then the statement of lemma is trivial.

\textbf{2.} To prove the equality $\pi_1=\pi_2$, it is enough to show that any
permutation $\tau\in \Symmgr_n$ is uniquely determined by $T_o^{\tau}$, $T_{cl}^{\tau}$,
and $T$.

For any $1\leq i\leq n$ there exists an $a\stackrel{t}{\in} T$,
satisfying $a'=q_2$ and $'a=i$. There is a unique maximal path $c$
from $T^{\tau}_{o}\sqcup T^{\tau}_{cl}$ that contains $a^{\tau}$
or $(a^t)^{\tau}$. Hence $c=a_1^{\tau}\cdots a_k^{\tau}\cdots
a_l^{\tau}$, where $a\stackrel{t}{=}a_k$ for some $k$, and
$a_1,\ldots,a_l\stackrel{t}{\in}T$.

Let $a=a_k$. If $k<l$, then $(a_{k+1}^{\tau})''=q_1=a_{k+1}''$ by the definition, which
implies $''a_{k+1}^{\tau}={}''a_{k+1}$. Moreover,
$\tau(i)={}'a^{\tau}={}''a_{k+1}^{\tau}={}''a_{k+1}$ and $''a_{k+1}$ is uniquely
determined by $T$. Therefore in the given case $\tau(i)$ is uniquely determined by
$T_o^{\tau}$, $T_{cl}^{\tau}$, and $T$. If $k=l$, then the path $c$ is open. Taking $a_1$
instead of $a_{k+1}$ and acting in the same manner as in the previous case, we obtain
that $\tau(i)={}''a_{1}$.

The case of $a=a_k^t$ can be treated analogously.
\end{proof}

\begin{lemma}\label{lemma_a_to_b}
Let $a,b\stackrel{t}{\in}T$ and $a\stackrel{t}{\neq}b$.

\begin{enumerate}
\item[a)] Suppose $a'=b'\in\{q_1,q_2\}$, $a''=b''\in\{q_1,q_2\}$.
Then for $\pi=(''a,{}''b)\circ ('a,{}'b)\in \Symmgr_n$ we have
$$T_{o}^{\pi}\stackrel{t}{\sim}T_{o}|_{a\to b,\;b\to a},\;%
T_{cl}^{\pi}\stackrel{ct}{\sim}T_{cl}|_{a\to b,\;b\to a}.$$

\item[b)] Suppose $a'=b'\in\{q_1,q_2\}$,
$a''=b''\not\in\{q_1,q_2\}$. Then for the transposition
$\pi=('a,{}'b)\in \Symmgr_n$ we have
$$T_{o}^{\pi}\stackrel{t}{\sim}T_{o}|_{a\to b,\;b\to a},\;%
T_{cl}^{\pi}\stackrel{ct}{\sim}T_{cl}.$$ %
\end{enumerate}
\end{lemma}
\begin{proof}{\bf a)}
Let $a'=b'=q_1$, $a''=b''=q_2$. Moreover, assume that
$'a,{}''a,{}'b,{}''b$ are pairwise different. Then $a,b$ belong to
some maximal paths $a_0$, $b_0$ in $T$, i.e., $a_0=a_1\cdots a_2 a
a_3\cdots a_4$ and $b_0=b_1\cdots b_2 b b_3\cdots b_4$, where
$a_i,b_j\stackrel{t}{\in}T$. It is not difficult to see that
maximal paths in $T^{\pi}$, containing $a^{\pi}$ or $b^{\pi}$, are
$a_1^{\pi}\cdots a_2^{\pi} b^{\pi} a_3^{\pi}\cdots a_4^{\pi}$ and
$b_1^{\pi}\cdots b_2^{\pi} a^{\pi} b_3^{\pi}\cdots b_4^{\pi}$. The
remaining maximal paths in $T^{\pi}$ are equivalent to paths in
$T$. The claim follows.

Other possibilities for $a,b$  can be treated in the similar
fashion.

{\bf b)} Analogous to part~a).
\end{proof}

Recall that Lemma~\ref{lemma_decomp} states that
$\F^0_T(X_1,\ldots,X_s)=\sum_{\tau\in \Symmgr_n}g(T,\tau)$, where
\begin{eq}\label{eq_g_T_tau}
g(T,\tau)=\sign(\tau)\F^0_{\widetilde{T}^{\tau}}(X_{\widetilde{T}^{\tau}})
\tr(T_{cl}^{\tau}).
\end{eq}

\begin{lemma}\label{lemma_g}
Let $\tau\in \Symmgr_n$, $a,b\stackrel{t}{\in} T^{\tau}$, $a\neq
b$, and $\ovphi{a}=\ovphi{b}$. If $a'=b'\in\{q_1,q_2\}$ and
$a''=b''\not\in \{q_1,q_2\}$, then %
$g(T, ('a,{}'b)\circ\tau)=g(T,\tau)$.
\end{lemma}
\begin{proof}
{\bf a)} Let $a'=b'=q_1$. For simplicity assume $\tau=id$ is the
identical permutation. Consider paths $u=u_1\cdots u_l$, $v$ in
$T$ such that $au,bv\stackrel{t}{\in}T_o$, $u_1,\ldots,
u_l\stackrel{t}{\in} T$ and set $k=a''=b''$ and $\xi=('a,{}'b)\in
\Symmgr_n$. There are two cases.

{\bf Case} $u_l\neq b^t$. We have
$$T_{o}^{\xi}\stackrel{t}{\sim} T_o\bigcup\{a^{\xi}v^{\xi},b^{\xi}u^{\xi}\}
\stackrel{t}{\backslash}
\{au,bv\},\;\; %
T_{cl}^{\xi}\stackrel{ct}{\sim}T_{cl}.$$ %
Therefore $$g(T,\xi)=-\sum\limits_{\pi_1\in
\Symmgr_{n_1},\ldots,\pi_{m-2}\in \Symmgr_{n_{m-2}}}
\sign(\pi_1)\cdots\sign(\pi_{m-2})\,
\tr(T_{cl})\ast$$ %
$$\prod_{d\in \widetilde{T},\;d\stackrel{t}{\neq}\widetilde{au},\widetilde{bv}} %
(X_{\ovphi{d}})_{\pi_{d''}(''d),\pi_{d'}('d)} %
\;\cdot\;
(X_{\ovphi{av}})_{\pi_{k}(''a),\pi_{v'}('v)}%
\;(X_{\ovphi{bu}})_{\pi_{k}(''b),\pi_{u'}('u)}.
$$
Observe that $X_{\ovphi{av}}=X_{\ovphi{bv}}$ and %
$X_{\ovphi{bu}}=X_{\ovphi{au}}$ and substitute $\pi_k\circ
(''a,{}''b)$ for $\pi_k$. Then $\pi_k(''a)$ turns into
$\pi_k(''b)$, $\pi_k(''b)$ turns into $\pi_k(''a)$, and the rest
of $\pi_p(q)$ does not change for $1\leq p\leq m-2$, $q\in
\Symmgr_{n_p}$. This proves the claim.

{\bf Case} $u_l= b^t$. Denote $u_1\cdots u_{l-1}$ by $w$. Then
$awb^t\stackrel{t}{\in} T_o$ and
$$T_{o}^{\xi}\stackrel{t}{\sim} T_o\bigcup\{a^{\xi}(w^{\xi})^t (b^{\xi})^t\}
\stackrel{t}{\backslash}
\{awb^t\},\;\; %
T_{cl}^{\xi}\stackrel{ct}{\sim}T_{cl}.$$ %
Therefore $$g(T,\xi)=-\sum\limits_{\pi_1\in
\Symmgr_{n_1},\ldots,\pi_{m-2}\in \Symmgr_{n_{m-2}}}
\sign(\pi_1)\cdots\sign(\pi_{m-2})
\tr(T_{cl})\ast$$ %
$$\prod_{d\in \widetilde{T},\;d\stackrel{t}{\neq}\widetilde{awb^t}}
(X_{\ovphi{d}})_{\pi_{d''}(''d),\pi_{d'}('d)} %
\;\cdot\;(X_{\ovphi{aw^t b^t}})_{\pi_{k}(''a),\pi_{k}(''b)}.$$%
Observe that $X_{\ovphi{aw^tb^t}}=X_{\ovphi{awb^t}}^t$.
Substitution of $\pi_k\circ (''a,{}''b)$ for $\pi_k$ completes the
proof.

The proof for an arbitrary $\tau$ is exactly the same as for
$\tau=id$.

{\bf b)} The proof in the case $a'=b'=q_2$ is analogous.
\end{proof}

\begin{lemma}\label{lemma_coef}
Suppose $\tau\in \Symmgr_n$. Then
\begin{enumerate}
\item[a)] $\#\Symmgr_{\tau}=\frac{
c^{(q_1,q_2)}_T}{\un{\be}!\,c(\un{\la}(\tau))}$, where $\tau\in
\Symmgr_{\un{\be},\un{\ga},\un{b},\un{c}}$,
$(\un{\be},\un{\ga},\un{b},\un{c})\in Q_T$.

\item[b)] If $\pi\in \Symmgr_{\tau}$, then $g(T,\pi)=g(T,\tau)$.
\end{enumerate}
\end{lemma}
\begin{proof}{\bf a)}
Given $\pi_1,\pi_2\in \Symmgr_n$, we write
$(T_{o}^{\pi_1},T_{cl}^{\pi_1})\stackrel{ct}{\sim}(T_{o}^{\pi_2},T_{cl}^{\pi_2})$
and say that these pairs are equivalent
if $T_{o}^{\pi_1}\stackrel{t}{\sim}T_{o}^{\pi_2}$ and %
$T_{cl}^{\pi_1}\stackrel{ct}{\sim}T_{cl}^{\pi_2}$. As in
Section~\ref{section_decomp_F0} define
$$A=\{ a\in T^{\tau}\,|\,a'>m-2\text{ or }a''>m-2\}.$$
Let $S$ be the set consisting of permutations $\pi\in \Symmgr_n$
for which there is a bijection $\La:A\to A$ satisfying
\begin{eq}\label{eq_item1}
\ovphi{\La(a)}=\ovphi{a} {\rm\;for\;all\;}a\in A, %
\end{eq}
\begin{eq}\label{eq_item2}%
(\ov{\La}(T_{o}^{\tau}),\ov{\La}(T_{cl}^{\tau}))\stackrel{ct}{\sim}%
(T_{o}^{\pi},T_{cl}^{\pi}), %
\end{eq}%
where $\ov{\La}(a)=\La(a)$, $\ov{\La}(a^t)=\La(a)^t$ for an $a\in
A$, and $\ov{\La}(a_1\cdots a_l)=\ov{\La}(a_1)\cdots
\ov{\La}(a_l)$ for $a_1,\ldots,a_l\stackrel{t}{\in} A$.

Definitions of $S$ and $\Symmgr_{\tau}$ imply the inclusion
$S\subset \Symmgr_{\tau}$. On the other hand, it is not difficult
to see that $\Symmgr_{\tau}\subset S$. Thus $\Symmgr_{\tau}=S$.

A bijection $\La:A\to A$ can be considered as a permutation of $A$
and can be written as a composition of transpositions on $A$. By
Lemma~\ref{lemma_a_to_b}, for every bijection $\La:A\to A$
satisfying~\Ref{eq_item1} there exists a $\pi\in \Symmgr_n$ such
that~\Ref{eq_item2} holds. This and Lemma~\ref{lemma_remark} show
that $\#S$ equals the number of $\stackrel{ct}{\sim}$-equivalence
classes of pairs
$(\ov{\La}(T_{o}^{\tau}),\ov{\La}(T_{cl}^{\tau}))$ for bijections
$\La:A\to A$ satisfying~\Ref{eq_item1}. Combinatorial argument
completes the proof.

{\bf b)} Using the notation of part~a), assume $\pi\in
\Symmgr_{\tau}$. Since $\Symmgr_{\tau}=S$, there exists a
bijection $\La:A\to A$ satisfying~\Ref{eq_item1}
and~\Ref{eq_item2}. As in part~a), for
such $\La$ we construct $\pi_0\in \Symmgr_n$ such that %
$$(\ov{\La}(T_{o}^{\tau}),\ov{\La}(T_{cl}^{\tau}))\stackrel{ct}{\sim}%
(T_{o}^{\pi_0\circ \tau},T_{cl}^{\pi_0\circ \tau}) $$ %
and $\pi_0$ is a composition of permutations of types
\begin{enumerate}
\item[1.] $(''a,{}''b)\circ('a,'b)$,

\item[2.] $('a,{}'b)$,
\end{enumerate}
where $a,b\stackrel{t}{\in} T^{\tau}$, $\ovphi{a}=\ovphi{b}$ (apply parts~a), b) of
Lemma~\ref{lemma_a_to_b} to $T^{\tau}$). Using Lemma~\ref{lemma_remark} we infer
$\pi=\pi_0\circ \tau$. For permutations of the first type we obtain $g(T,
(''a,{}''b)\circ('a,'b)\circ\tau)=g(T,\tau)$; for permutations of the second type we use
Lemma~\ref{lemma_g}.
\end{proof}

Now we can express $\F_T(X_1,\ldots,X_s)$ in a suitable form. By
Lemma~\ref{lemma_decomp}, %
$\F_T(X_1,\ldots,X_s)=\frac{1}{c_T}\sum_{\tau\in
\Symmgr_n}g(T,\tau)$. %
Using Lemma~\ref{lemma_abcd} and notation from the beginning of
this section, we obtain
$$\F_T(X_1,\ldots,X_s)=\frac{1}{c_T}\sum_{(\un{\be},\un{\ga},\un{b},\un{c})\in Q_T}\;%
\sum_{\pi\in\ov{\Symmgr}_{\un{\be},\un{\ga},\un{b},\un{c}}}\;\sum_{\tau\in
\Symmgr_{\pi}}g(T,\tau).$$ %
Lemma~\ref{lemma_coef} yields %
$$\F_T(X_1,\ldots,X_s)$$
$$=\frac{1}{c^{(q_1,q_2)}_T c^{(rest)}_T}\sum_{(\un{\be},\un{\ga},\un{b},\un{c})\in Q_T}%
\sum_{\pi\in\ov{\Symmgr}_{\un{\be},\un{\ga},\un{b},\un{c}}}
\frac{c_T^{(q_1,q_2)}}{\un{\be}!\,c(\un{\la}(\pi))}
\sign(\pi)\F^0_{\widetilde{T}^{\pi}}(X_{\widetilde{T}^{\pi}})
\tr(T_{cl}^{\pi}).$$ %
Let us recall that
$c_{\widetilde{T}^{\pi}}=c_T^{(rest)}c^{(\pi)}_T$. Since
$c^{(\pi)}_T=\un{\be}!$ for
$\pi\in\ov{\Symmgr}_{\un{\be},\un{\ga},\un{b},\un{c}}$ we conclude %
\begin{eq}\label{eq_coef}%
\F_T(X_1,\ldots,X_s)=\sum_{(\un{\be},\un{\ga},\un{b},\un{c})\in Q_T}\;%
\sum_{\pi\in\ov{\Symmgr}_{\un{\be},\un{\ga},\un{b},\un{c}}}
\sign(\pi)\frac{1}{c(\un{\la}(\pi))}\tr(T_{cl}^{\pi})
\F_{\widetilde{T}^{\pi}}(X_{\widetilde{T}^{\pi}}).%
\end{eq}

\section{Traces and $\si_k$}\label{section_tr}
Continue working under the same assumptions as in the previous
section.

For $\ga\in\NN$, $\un{\nu}=(\nu_1,\ldots,\nu_l)\vdash \ga$ denote
by $\rho_{\un{\nu}}\in \Symmgr_{\ga}$ the composition of cycles
$$\rho_{\un{\nu}}=(1,2,\ldots,\nu_1)\circ (\nu_1+1,\ldots,\nu_1+\nu_2)\circ\cdots\circ
(\nu_1+\cdots+\nu_{l-1},\ldots,\ga).$$ %
Note that $\sign(\rho_{\un{\nu}})=(-1)^{\ga-\#\un{\nu}}$.

Let $(\un{\be},\un{\ga},\un{b},\un{c})\in Q_T$,
$\xi=\xi_{\un{\be},\un{\ga},\un{b},\un{c}}$, and $\#\un{\ga}=q$.
Without loss of
generality we can assume that %
$$T_{cl}^{\xi}=\Lus d_{11},\ldots,d_{1\ga_1},\ldots,d_{q1},\ldots,d_{q\ga_q}\Rus,
\quad{\rm where}$$ %
\begin{eq}\label{eq_d}
d_{ij}\in\M_T,\;\;\ \ovphi{d_{ij}}=c_i,\;\; %
d_{ij}''=q_1,\; d_{ij}'=q_2,\;{}''d_{ij}={}'d_{ij}
\end{eq} %
for any $1\leq i\leq q$ and $1\leq j\leq \ga_i$. There is a unique
bijection $\phi_i:[1,\ga_i]\to \{{}'d_{i1},\ldots,{}'d_{i\ga_i}\}$
that is a monotone increasing map, i.e., $\phi_i(j_1)<\phi_i(j_2)$ for $j_1<j_2$. %

Let $S=\{\pi_{1,\un{\nu}_1}\circ \cdots\circ
\pi_{q,\un{\nu}_q}\circ \xi\;|\;
\un{\nu}_1\vdash\ga_1,\ldots,\un{\nu}_q\vdash\ga_q\}$, where
the permutation $\pi_{i,\un{\nu}_i}\in \Symmgr_n$ is defined by %
$$\pi_{i,\un{\nu}_i}(k)=\left\{ %
\begin{array}{cl}
\phi_i\circ \rho_{\un{\nu}_i}\circ \phi^{-1}_i(k) ,&\text{if }
k\in\image{\phi_i} \\
k,&\text{otherwise}
\end{array}
\right..$$%
Then
\begin{enumerate}
\item[$\bullet$] By~\Ref{eq_la}, $\un{\la}(\pi_{1,\un{\nu}_1}\circ
\cdots\circ \pi_{q,\un{\nu}_q}\circ
\xi)=(\un{\nu}_1,\ldots,\un{\nu}_q)$ is a multi-partition of
$\un{\ga}$.

\item[$\bullet$] Suppose $\pi,\tau\in
\Symmgr_{\un{\be},\un{\ga},\un{b},\un{c}}$. Then
$\Symmgr_{\tau}=\Symmgr_{\pi}$ if and only if
$\un{\la}(\tau)=\un{\la}(\pi)$.
\end{enumerate}
These remarks imply that
$\Symmgr_{\un{\be},\un{\ga},\un{b},\un{c}}=\bigsqcup_{\tau\in S}
\Symmgr_{\tau}$, hence without loss of generality we can assume
that $\ov{\Symmgr}_{\un{\be},\un{\ga},\un{b},\un{c}}=S$
(see~\Ref{eq_ovS}).

Note that
\begin{enumerate}
\item[$\bullet$] Since
$\sign(\pi_{i,\un{\nu}_i})=\sign(\rho_{\un{\nu}_i})$, we have the
equality $\sign(\pi_{1,\un{\nu}_1}\circ \cdots\circ
\pi_{q,\un{\nu}_q}\circ \xi)=\sign(\rho_{{\un{\nu}}_1}) \cdots
\sign(\rho_{{\un{\nu}}_q})\sign(\xi)$.

\item[$\bullet$] $\tr(T_{cl}^{\pi_{1,\un{\nu}_1}\circ \cdots\circ
\pi_{q,\un{\nu}_q}\circ \xi})=\prod\limits_{i=1}^q
\prod\limits_{j=1}^{\#\un{\nu}_i} \tr(X_{c_i}^{\nu_{ij}})$. %

\item[$\bullet$] Since $T_{o}^{\eta}\stackrel{t}{=}T_{o}^{\xi}$
for $\eta= {\pi_{1,\un{\nu}_1}\circ \cdots\circ
\pi_{q,\un{\nu}_q}\circ \xi}$, it is easy to see that
$(\widetilde{T}^{\eta},X_{\widetilde{T}^{\eta}})=
(\widetilde{T}^{\xi},X_{\widetilde{T}^{\xi}})$.

\end{enumerate}
Therefore equality~\Ref{eq_coef} yields $\F_T(X_1,\ldots,X_s)=$ %
\begin{eq}\label{eq_F}
\sum_{(\un{\be},\un{\ga},\un{b},\un{c})\in Q_T}%
\sign(\xi_{\un{\be},\un{\ga},\un{b},\un{c}})
\F_{\widetilde{T}^{\xi_{\un{\be},\un{\ga},\un{b},\un{c}}}}
(X_{\widetilde{T}^{\xi_{\un{\be},\un{\ga},\un{b},\un{c}}}})
\prod_{i=1}^{\#\un{\ga}}\; %
\left( \sum_{\un{\nu}\vdash \ga_i}
(-1)^{\ga_i-\#\un{\nu}}\frac{1}{c(\un{\nu})}\;\prod_{j=1}^{\#\un{\nu}}
\tr(X_{c_i}^{\nu_j}) \right).
\end{eq} %
This formula together with the following well known lemma completes the proof of
Theorem~\ref{theo_decomp}.

\begin{lemma}\label{lemma_zeta}
Assume $K=\QQ$. If $k,n\in \NN$, $0<k\leq n$, and $c$ is an
$n\times n$ matrix, then
$$\si_{k}(c)=\sum_{\un{\nu}\vdash k}
(-1)^{k-\#\un{\nu}}\,\frac{1}{c(\un{\nu})}\;\prod_{j=1}^{\#\un{\nu}}\tr(c^{\nu_j}),$$ %
where $c(\un{\nu})$ is given by~\Ref{eq_c_la}.
\end{lemma}
\bigskip

Applying~\Ref{eq_F} for the tableau $T$ from part~$3$ of Example~\ref{ex1} we obtain an
alternative proof of this lemma.

%

\section{Corollaries}\label{section_examples}


The decomposition formula generalizes Amitsur's formula for the
determinant (see~\cite{Amitsur}).

\begin{cor}\label{cor_Amitsur}
{\bf(Amitsur).} For $k_1+\cdots+k_s=n$ and $n\times n$ matrices
$X_1,\ldots,X_s$
we have %
$$\det\nolimits_{k_1,\ldots,k_s}(X_1,\ldots,X_s)=\sum
(-1)^{n+(i_1+\cdots+i_q)}
\si_{i_1}(c_1)\cdots\si_{i_q}(c_q),$$ %
where the sum ranges over all $\stackrel{c}{=}$-equivalent classes
of multisets  \\%
$\Lus\underbrace{c_1,\ldots,c_1}_{i_1},\ldots,\underbrace{c_q,\ldots,c_q}_{i_q}\Rus$
such that
\begin{enumerate}
\item[$\bullet$] $c_1,\ldots,c_q$ are primitive words in letters
$X_1,\ldots,X_s$ such that $c_1,\ldots,c_q$ are pairwise different
with respect to $\stackrel{c}{=}$;

\item[$\bullet$] $i_1\deg_{X_j}c_1+\cdots+i_q\deg_{X_j}c_q=n_j$
for any $1\leq j\leq s$.
\end{enumerate}
\end{cor}
\begin{proof}
Consider the tableau with substitution $(T,(X_1,\ldots,X_s))$ of dimension $(n,n)$ from
part~2 of Example~\ref{ex1} and assume $\#\{a\in T\,|\, \ovphi{a}=j\}=k_j$ for any $1\leq
j\leq s$. Then $\F_T(X_1,\ldots,X_s)=\det_{k_1,\ldots,k_s}(X_1,\ldots,X_s)$ and the claim
follows easily from the decomposition formula.
\end{proof}

\begin{cor}\label{cor_PP}
Let $Y$, $Z$ be $n\times n$ matrices and $n$ is even. Then
$\P(Y)\P(Z)=$ %
$$\sum (-1)^{i_1(\deg_Y c_1+\deg_Z c_1 +1)+\cdots+i_q(\deg_Y c_q+\deg_Z c_q +1)}
\si_{i_1}(c_1)\cdots\si_{i_q}(c_q),$$ %
where the sum ranges over $\stackrel{ct}{=}$-equivalent
classes of multisets  \\ %
$\Lus\underbrace{c_1,\ldots,c_1}_{i_1},\ldots,\underbrace{c_q,\ldots,c_q}_{i_q}\Rus$
such that
\begin{enumerate}
\item[$\bullet$] $c_1,\ldots,c_q$ are primitive words in letters
$Y,Y^t,Z,Z^t$ such that $c_1,\ldots,c_q$ are pairwise different
with respect to $\stackrel{ct}{=}$;

\item[$\bullet$] words $c_1,\ldots,c_q$ are products of the words
$AB$, where $A\in\{Y,Y^t\}$ and $B\in\{Z,Z^t\}$.

\item[$\bullet$]$i_1(\deg_Y c_1+\deg_{Y^t} c_1)+\cdots+%
i_q(\deg_Y c_q+\deg_{Y^t} c_q) =n/2$, \\
$i_1(\deg_Z c_1+\deg_{Z^t} c_1)+\cdots+%
i_q(\deg_Z c_q+\deg_{Z^t}c_q) =n/2$.
\end{enumerate}
\end{cor}
\begin{proof}
Consider the tableau with substitution $(T,(Y,Z))$ of dimension $(n,n)$ from part~4 of
Example~\ref{ex1}, where $t=0$, $r=s=n/2$. Then $\F_T(Y,Z)=\P(Y)\P(Z)$ and the formula
follows from Theorem~\ref{theo_decomp}.
\end{proof}

\example\label{ex4} 1. Let $Y$, $Z$ be $2\times 2$ matrices. Then
$$\P(Y)\P(Z)=\tr(YZ^t) - \tr(YZ).$$

2. Suppose $Y$, $Z$ are $4\times 4$ matrices. Then
$\P(Y)\P(Z)=$ %
$$\si_2(YZ)+\si_2(YZ^t)-%
\tr(YZ)\tr(YZ^t)+\tr(YZYZ^t)+\tr(YZY^tZ)-\tr(YZY^tZ^t).$$

\begin{remark}
In the same manner as in Corollary~\ref{cor_PP} we can also write a formula for a product
of two partial linearizations of pfaffians.
\end{remark}


\bigskip
\noindent{\bf Acknowledgements.} The first version of this paper was written during
author's visit to University of Antwerp, sponsored by Marie Curie Research Training
Network Liegrits. The author is grateful for this support. The author would like to thank
Fred Van Oystaeyen for his hospitality. This research was also supported by RFFI
05-01-00057.


\end{document}